\documentclass[a4paper,12pt]{article}
\usepackage{amsmath}
\numberwithin{equation}{section}
\newtheorem{theorem}{Theorem}[section]
\newtheorem{proposition}{Proposition}[section]
\newtheorem{corollary}{Corollary}[section]
\newtheorem{lemma}{Lemma}[section]
\newtheorem{definition}{Definition}[section]
\title{\textbf{Deformed Infinite series metric in Cartan Spaces}}
\date{}
\begin{document}
\maketitle
\begin{center}
\textbf{\author{Brijesh Kumar Tripathi and $ ^{*} $V. K. Chaubey}}
\end{center}

\begin{center}
	Department of Mathematics, L.D. College of Engineering, Ahmedabad \\* (Gujarat)-380015, India, E-mail:brijeshkumartripathi4@gmail.com
\end{center}
\begin{center}
$ ^{*} $ Department of Applied Sciences, Buddha Institute of Technology, GIDA \\* Gorakhpur, U.P., 273209, India, E-Mail: vkcoct@gmail.com
\end{center}
\begin{abstract}
Igarashi introduce the concept of $(\alpha, \beta)$-metric in Cartan space $\ell^{n}$ analogously to one in Finsler space and obtained the basic important geometric properties and also investigate the special class of the space with $(\alpha, \beta)$-metric in $\ell^{n}$ in terms of $ 'invariants' $. In the present paper we determine the $ 'invariants' $ in two different cases of deformed infinite series metric which characterize the special classes of Cartan spaces $\ell^{n}$. 
\end{abstract}
\section{Introduction}
\hspace{10pt} In 2004 Lee and Park \cite{LP} introduced the concept of r-th series $(\alpha,\beta)$-metric where r is varies from $ 0, 1, 2, ..., \infty $ and give very interesting example of special $(\alpha, \beta)$-metric for the different values of r such as one-form metric, Randers metric, combination of Kropina and Randers metric, infinite series metric etc.\\

In 1994, Igarashi \cite{Iga1, Iga2} introduce the concept of $(\alpha, \beta)$-metric in Cartan space $\ell^{n}$ analogously to one in Finsler space and obtained the basic important geometric properties and also investigate the special class of the space with $(\alpha, \beta)$-metric in $\ell^{n}$ in terms of $ 'invariants' $. The classes which he obtained includes the spaces corresponding to Randers and Kropina space. Further he characterizes these spacial classes by means of $ 'invariants' $ in case of Finsler theory.\\

In the present paper we determine the $ 'invariants' $ in two different cases of deformed infinite series metric in which firest metric is defined as the product of infinite series and Riemannian metric another one is the product of infinite series and one-form metric. Further we characterize the special classes of Cartan spaces $\ell^{n}$ in case of these two metrics and also investigate the relation under which "invariants" are characterized as the special classes of $\ell^{n}$

\section{Preliminaries}
E. Cartan \cite{Car} introduced the concept of a Cartan Space, where the measure of its hypersurface element $(x,y)$ is given a priori by homogeneous function $F(x,y)$ of degree one in y,i.e.,the "area" of a domain on hypersurface $S_{n-1} : x^{i} = x^{i}(v^{1},v^{2},v^{3},.....,v^{n}), \;\;\; i=1,2,3,....,n$ is given by 
\begin{equation}
S=\underbrace{\idotsint}F(x,y)dv^{1}dv^{2}dv^{3}.....dv^{n-1}
\end{equation}
where $ y = (y_{i})$ is the determinant of $ (n-1, n-1)$ minor matrix omitted ith row of $(n, n-1)$ matrix $(\frac{\partial x^{i}}{\partial v^{\alpha}})$, \;\;\; $\alpha=1,2,3.....n-1$. In this space, we obtain the fundamental tensor by
\begin{equation}
g^{ij}=G^{-\frac{1}{n-1}},\;\;\; G=det\lVert G^{ij}\rVert, \;\;\; G^{ij}=\frac{\partial^{2}(\frac{1}{2}F^{2})}{\partial y_{i}\partial y_{j}}, \;\;\; (y_{i})\neq 0.
\end{equation}
As the special case for the fundamental tensor $a_{ij}$ of Riemannian space, we can find the $(n-1)$ dimensional area of a domain on hypersurface such that
\begin{center}
$S=\underbrace{\idotsint}\sqrt{det\lVert a_{ij}(x)\frac{\partial x^{i}}{\partial v^{\alpha}}\frac{\partial x^{j}}{\partial v^{\beta}}\rVert}dv^{1}dv^{2}dv^{3}.....dv^{n-1}$
\end{center}
hence it is clear that Riemannian space is a special case of Cartan space. On the other hand, Cartan space is considered the dual notion of Finsler Space. Further the relation between both spaces is studied by L. Berwald \cite{Be} in early days, afterwards, by H. Rund \cite{Ra} and F.Brickel \cite{Br}. Recently R. Miron \cite{Mi1, Mi2} established new Carton geometry which shows totally different feature in the form of particularization the Hamilton space which defined as:
\begin{definition}
A Cartan space is a Hamilton space $\mathcal{H}^{n}=\{M,\underbrace{H(x,y)}\}$ in which the fundamental function $ H(x,y)$ is positively 2-homogeneous in $ y_{i}$ on $T^{*}M$. We denote it by $\ell^{n}$.
\end{definition}
The fundamental tensor field of $\ell^{n}$ and its reciprocal $g_{ij}(x,y)$ is given by
\begin{equation}
g^{ij}(x,y)=\frac{1}{2}\dot{\partial^{i}}\dot{\partial^{j}}H,
\end{equation}
\begin{equation}
g_{ij}(x,y)g^{jk}(x,y)=\delta_{i}^{k}
\end{equation}
The homogeneity of $H(x,y)$ is expressed by
\begin{equation}
y_{j}\dot{\partial^{j}}H=2H,\,\,  which\,\, also\,\, implies\,\, H=g^{ij}y_{i}y_{j}
\end{equation}
where $g^{ij}(x,y)$ and its reciprocal $g_{ij}(x,y)$ are both symmetric and homogeneous of degree 0 in $y_{i}$.\\\\
On the other hand, the Finsler spaces with $(\alpha,\beta)$-metric were considered by G. Randers \cite{Ra}, V. K. Kropina \cite{Kro} and M. Matsumoto \cite{M1, M2, M3}, especially the last paper shows the great success for investigation of these spaces.\\

In \cite{Mi1}, R. Miron expected the existence of Randers type metric:
\begin{equation}
H(x,y)=\lbrace\alpha(x,y)+\beta(x,y)\rbrace^{2}
\end{equation}
and of Kropina's one :
\begin{equation}
H(x,y)=\{\frac{[\alpha(x,y)]^{2}}{\beta(x,y)}\}^{2},\,\,\,\,\,\,\,\,\beta(x,y)\neq 0.
\end{equation}
in Cartan spaces. Here we have put as
\begin{equation}
\alpha^{2}(x,y)=a^{ij}(x) y_{i} y_{j},\,\,\,\,\beta(x,y)=b^{i}(x)y_{i},
\end{equation}
$a^{ij}(x)$ being a Riemannian metric on the base manifold M and $b^{i}(x)$ a vector field on M such that $\beta>0$ on a region of $T^{*}M\stackrel{def}{=} T^{*}M-\lbrace 0 \rbrace $.

\section{Cartan spaces with $(\alpha,\beta)$-metric.}
Cartan spaces with $(\alpha,\beta)$-metric \cite{Iga1, Iga2, Mi1} can be defiend as
\begin{definition}
A Cartan space $\ell^{n}=\{M, H(x,y)\}$  is known as Cartan space with $(\alpha, \beta)$-metric if its fundamental metric $H(x,y)$ is a function of $\alpha(x,y)$ and $\beta(x,y)$ only i.e.
\begin{equation}
H(x,y)=\breve{H}\{\alpha(x,y), \beta(x,y)\}
\end{equation}
\end{definition}
It is clear that $\breve{H}$ satisfy the conditions imposed to the function $H(x,y)$ as a fundamental function for $\ell^{n}$. Then

\begin{equation}
t\alpha(x,y)=\alpha(x,ty),\;\; t\beta(x,y)=\beta(x,ty),\;\; H(x,ty)=t^{2}H(x,y)\;\; t>0.
\end{equation}
It follows:
\begin{proposition}
The function $\breve{H}(\alpha(x,y),\beta(x,y))$ is positively homogeneous of degree 2 in both $\alpha$ and $\beta$.
\end{proposition}
By this reason, there maybe no confusion if we adopt the notation $H(\alpha,\beta)$ itself instead of $\breve{H}(\alpha,\beta)$. Also we write such that
\begin{equation}
H_{\alpha}=\frac{\partial H}{\partial \alpha},\,\,H_{\beta}=\frac{\partial H}{\partial \beta},\,\,H_{\alpha\beta}=\frac{\partial^{2} H}{\partial \alpha\partial \beta},\,\,etc
\end{equation}

\begin{proposition}
The following identities hold :
\begin{equation}
\alpha H_{\alpha\alpha}+\beta H_{\beta\beta}=H,\,\,\alpha H_{\alpha\beta}+\beta H_{\beta\beta}=H_{\beta},
\end{equation}
\begin{equation}
 \alpha^{2}H_{\alpha \alpha}+2\alpha\beta H_{\alpha\beta}+\beta^{2} H_{\beta\beta}=\alpha H_{\alpha}+\beta H_{\beta}=2H   \nonumber .
\end{equation}
\end{proposition}

Differentiating $ \alpha $ and $ \beta $ with respect to $y_{i}$ we have
\begin{equation}
\dot{\partial^{i}} \alpha=\alpha^{-1}a^{ij}y_{j}=\alpha^{-1}Y^{i},\;\;\;\;\;\dot{\partial^{i}}\beta=b^{i}(x),
\end{equation}
where
\begin{equation} 
Y^i(x,y)=a^{ij}(x)y_{i},\;\;\; \{Y=(Y^{i})\}\ne 0 
\end{equation}
and the vector field $ Y^{i}$ satisfies the relation
\begin{equation} 
 Y^{i}y_{i}=\alpha^{2},
\end{equation}
Further let
\begin{equation} 
B_{i}(x)=a_{ij}(x)b^{j}(x),\;\;\;B^{2}(x)=a^{ij}B_{i}B_{j}=a_{ij}b^{i}b^{j},
\end{equation}
\begin{equation} 
y^{i}(x,y)=g^{ij}(x,y)y_{j}=\frac{1}{2}\frac{\partial H}{\partial y_{i}}
\end{equation}
Also, we have the relation similar to (3.7):
\begin{equation} 
y^{i}y_{i}=H(x,y).
\end{equation}
Differentiating (3.6) and (3.9) by $y_{i}$ succeedingly, we have 
\begin{equation} 
\dot{\partial ^{j}}Y^{i}=a^{ij}(x),\;\;\; \dot{\partial ^{k}}\dot{\partial^{j}}Y^{i}=0,
\end{equation} 
\begin{equation}
\dot{\partial^{j}}y^{i}=g^{ij},\;\;\; \dot{\partial^{k}}\dot{\partial^{j}}y^{i}=\dot{\partial^{k}}g^{ij}=-2C^{ijk},
\end{equation}
And using the same manner for (3.5), we get
\begin{equation} 
\dot{\partial^{j}}\dot{\partial^{i}}\alpha=\alpha^{-1} a^{ij}(x)-\alpha^{-3}Y^{i}Y^{j},\;\;\dot{\partial^{j}}\dot{\partial^{i}}\beta=0,\;\;\;\dot{\partial^{j}}\dot{\partial^{i}}(\frac{1}{2}\alpha^{2})=a^{ij}.
\end{equation}
On account of (3.5) and $y^{i}=\frac{1}{2}(H_{\alpha}\dot{\partial^{i}}\alpha+H_{\beta}\dot{\partial^{i}}\beta)$, we have
\begin{lemma}
The Liouville vector field $y^{i}$ is expressed in the form
\begin{equation}
y^{i}=\rho_{1}b^{i}+\rho Y^{i},
\end{equation}
where
\begin{equation}
\rho_{1}=\frac{1}{2}H_{\beta},\;\;\;\; \rho=\frac{1}{2\alpha}H_{\alpha}
\end{equation}
\end{lemma}
Taking into account that the relation $\dot{\partial^{i}}\rho =\frac{\partial\rho}{\partial \alpha}\dot{\partial^{i}}\alpha +\frac{\partial\rho}{\partial \beta}\dot{\partial^{i}}\beta $ holds, we have
\begin{lemma}
The quantities $\rho_{1}$ and $\rho$ satify the relations
\begin{equation}
\dot{\partial^{i}}\rho_{1} =\rho_{0}b^{i}+ \rho_{-1}Y^{i},\;\;\;\; \dot{\partial^{i}}\rho =\rho_{-1}b^{i}+ \rho_{-2}Y^{i}
\end{equation}
respectively, where 
\begin{equation}
\rho_{0}=\frac{1}{2}H_{\beta\beta},\;\;\; \rho_{-1}=\frac{1}{2\alpha}H_{\alpha\beta},\;\;\;\rho_{-2}=\frac{1}{2\alpha^{2}}(H_{\alpha\alpha}-\alpha^{-1}H_{\alpha})
\end{equation}
\end{lemma}
Contracting $\dot{\partial^{i}}\rho_{1}$,\;\;$ \dot{\partial^{i}}\rho$ in (3.16) by $y_{i}$,we have\\
\begin{equation} 
y_{i}\dot{\partial^{i}}\rho_1=\alpha^{2}\rho_{-1}+\beta\rho_{0}=\rho_{1},\;\;\;y_{i}\dot{\partial^{i}}\rho=\alpha^{2}\rho_{-2}+\beta\rho_{-1}=0.
\end{equation} 
Analogously to deduction of Lemma (3.2), we have also
\begin{lemma}
The quantities $\rho_{0}$ and $\rho_{-1}$ satisfy the relations
\begin{equation}
\dot{\partial^{i}}\rho_{0} =r_{-1}b^{i}+ r_{-2}Y^{i},\;\;\;\; \dot{\partial^{i}}\rho_{-1} =r_{-2}b^{i}+ r_{-3}Y^{i},
\end{equation} 
respectively, where
\begin{equation}
r_{-1}=\frac{1}{2}H_{\beta\beta\beta},\;\;\; r_{-2}=\frac{1}{2\alpha}H_{\alpha\beta\beta},\;\;\;r_{-3}=\frac{1}{2\alpha^{2}}(H_{\alpha\alpha\beta}-\alpha^{-1}H_{\alpha\beta}).
\end{equation}
\end{lemma}
\;\;\; Corresponding to (3.18), it follows
\begin{equation}
y_{i}\dot{\partial^{i}}\rho_0=\alpha^{2}r_{-2}+\beta r_{-1}=0,\;\;\;y_{i}\dot{\partial^{i}}\rho_{-1}=\alpha^{2} r_{-3}+\beta r_{-2}=-\rho_{-1}.
\end{equation} 
Furthermore for $\rho_{-2}$, we have
\begin{lemma}
The quantities $\rho_{-2}$ satisfies the relations
\begin{equation}
\dot{\partial^{i}}\rho_{-2} =r_{-3}b^{i}+ r_{-4}Y^{i}
\end{equation}
\begin{equation}
r_{-4}=\frac{1}{2\alpha^{3}}(H_{\alpha\alpha\alpha}-3\alpha^{-1}H_{\alpha\alpha}+3\alpha^{-2}H_{\alpha}).
\end{equation}
\end{lemma} 
Also following homogeneity holds good:
\begin{equation}
y_{i}\dot{\partial^{i}}\rho_{-2}=\alpha^{2}r_{-4}+\beta r_{-3}=-2\rho_{-2}.
\end{equation}
It is easy to conclude for the scalars (or invariants) $ \rho_{1},\rho,\rho_{0},\rho_{-1},.......,r_{-1}, r{-2},......$ in the above lemmas that The subscript of $\rho$,s and $r$,s represent degree of their own homogeneity in $(\alpha,\beta)$ or $y_{i}$,where the $\rho$ without subscript means of degree 0.\\

We have these properties from expressions in (3.18),(3.19) and (3.21) and the following relations
\begin{equation}
y_{i}\dot{\partial^{i}}r_{-1}=-r_{-1},\;\;\;y_{i}\dot{\partial^{i}}r_{-2}=-2r_{-2},
\end{equation}
\begin{equation}
y_{i}\dot{\partial^{i}}r_{-3}=-3r_{-3},\;\;\;y_{i}\dot{\partial^{i}}r_{-4}=-4r_{-4},\nonumber
\end{equation}
because of the homogeneity of
\begin{equation}
\alpha H_{\alpha\alpha\alpha\alpha}+\beta H_{\alpha\alpha\alpha\beta}=-H{\alpha\alpha\alpha},\;\;\; \alpha H_{\alpha\alpha\alpha\beta}+\beta H_{\alpha\alpha\beta\beta}=-H{\alpha\alpha\beta},
\end{equation}
\begin{equation}
\;\;\alpha H_{\alpha\alpha\beta\beta}+\beta H_{\alpha\beta\beta\beta}=-H{\alpha\beta\beta},\;\;\; \alpha H_{\alpha\beta\beta\beta}+\beta H_{\beta\beta\beta\beta}=-H{\beta\beta\beta}.\nonumber
\end{equation}

We shall use the previous results for study the fundamental geometric objects of the space $\ell^{n}$ with $(\alpha,\beta)$-metric.All these scalar functions $ \rho_{1},\rho,\rho_{0},\rho_{-1},......$ as well as $r_{-1},r_{-2},.....$ will be called the invariants of the Cartan space $\ell^{n}$ with the fundamental function $H(\alpha,\beta)$.\\
 
\section{The fundamental tensor of the space $\ell^{n}$ with $(\alpha,\beta)$-metric.}
  
  We need to derive the fundamental tensor from the fundamental function $H(x,y)$ of the Cartan space $\ell^{n}$.
  \begin{theorem}
  The fundamental tensor $g^{ij}$ of Cartan space $\ell^{n}$ with $(\alpha,\beta)$-metric is given by
  \begin{equation}
  g^{ij}=\rho a^{ij}+\rho_{0}b^{i}b^{j}+\rho_{-1}(b^{i}Y^{j}+b^{j}Y^{i})+\rho_{-2}Y^{i}Y^{j},
  \end{equation}
  where $\rho,\rho_{0},\rho_{-1}\rho_{-2}$ are the invariants given by (2.15) and (2.17).
  \end{theorem}
\textbf{{Proof.}}\;\;\;
Making use of (3.12) and (3.14),we have
\begin{center}
$ g^{ij}=\dot{\partial^{j}}y^{i}=\dot{\partial^{j}}(\rho_{1}b^{i}+\rho Y^{i})=(\dot{\partial^{j}}\rho_{1})b^{i}+(\dot{\partial^{j}}\rho)Y^{i}+\rho a^{ij}$.
\end{center}
Taking into account Lemma (3.2), we have (3.1).\;\;\;\;\;\;\;\;\;\;Q.E.D.\\
In order to check the fitness of this tensor $g^{ij}$ for the fundamental tensor of $(\alpha,\beta)$-metric,we verify the homogeneity of $g^{ij}$.Contracting $g^{ij}$ by $y_{i}$ and $y_{j}$,we have
\begin{center}
$g^{ij}y_{i}y_{j}=\frac{1}{2}\lbrace \alpha^{-1}H_{\alpha}a^{ij}y_{i}y_{j}+(\alpha^{-2}H_{\alpha\alpha}-\alpha^{-3}H_{\alpha})Y^{i}Y^{j}y_{i}y_{j}$
$+\alpha^{-1}H_{\alpha\beta}(Y^{i}y_{i}b^{j}y_{j}+b^{i}y_{i}Y^{j}y_{j})+H_{\beta\beta}b^{i}y_{i}b^{j}y_{j}\rbrace$
$=\frac{1}{2}\lbrace \alpha H_{\alpha}+(\alpha^{-2}H_{\alpha\alpha}-\alpha^{-3}H_{\alpha})\alpha^{4}+\alpha^{-1}H_{\alpha\beta}.2\alpha^{2}\beta+H_{\beta\beta} \}$\\
$=\frac{1}{2}.2H=H.$
\end{center}
which shows that our conclusion is right.\\
Let us rewrite this expression in the form
\begin{equation}
g^{ij}=A^{ij}+C^{i}C^{j},
\end{equation}
where
\begin{equation}
A^{ij}=\rho a^{ij},\;\;\;C^{i}=q_{0}b^{i}+q^{-1}Y^{i},
\end{equation}
\begin{equation}
\rho_{0}=q^{2}_{0},\;\;\;\; \rho_{-1}=q_{0}q_{-1},\;\;\; \rho_{-2}=q_{-1}^{2}
\end{equation}
and
\begin{equation}
\rho_{0}\rho_{-2}=\rho_{-1}^{2}
\end{equation}
The reciprocal tensor $g_{ij}$ of $g^{ij}$ are given by
\begin{equation}
g_{ij}=A_{ij}-\frac{1}{1+C^{2}}C_{i}C_{j},
\end{equation}
where
\begin{equation}
det\lVert A_{ij}\rVert=(1+C^{2})det\lVert A^{ij}\rVert\;\;\;\;(if  1+C^{2}\ne 0),
\end{equation}
and $A_{ij}$,$C_{i}$ are given by
\begin{equation}
A_{ij}A^{jk}=\delta^{i}_{k},\;\;\; C^{i}C_{j}=\delta^{i}_{j}C^{2},\;\;\; C^{i}=A^{ij}C_{j},\;\;\;C_{i}=A_{ij}C^{j}.
\end{equation}
Consequently, we have $g^{ij}(x,y)g_{jk}(x,y)=\delta^{i}_{k}$,and\\
{(3.2')}\space\;\;\;\;\;\;\;\;\; $rank\lVert g^{ij}(x,y)\rVert=n$\\
because of\\
{(3.2'')\;\;\;\;\;\;$det\lVert g^{ij}(x,y)\rVert=(1+C^{2})det\lVert A^{ij}\rVert=(1+C^{2})det\lVert a^{ij}(x)\rVert\ne 0$\\
The following relations are useful afterwards :
\begin{equation}
A_{ij}=\frac{1}{\rho}a_{ij},\;\;\;det\lVert A^{ij}\rVert=\rho^{n}det\lVert a^{ij}\rVert ,\;\;(A=det\lVert a^{ij}\rVert\ne 0),
\end{equation}
\begin{equation}
C^{2}=\frac{1}{\rho}(\rho_{0}B^{2}+\rho_{-1}\beta),\;\;\;det\lVert g^{ij}\rVert=\rho^{n-1}\tau,
\end{equation}
where we use the notations
\begin{equation}
B^{2}=a^{ij}b_{i}b_{j},\;\;\; \tau=\rho+\rho_{0} B^{2}+\rho_{-1}\beta
\end{equation}
Therefore we can prove without difficulty:
\begin{proposition}
The covariant form of the fundamental tensor is given by
\begin{equation}
g_{ij}=\sigma a_{ij}-\sigma_{0} B_{i}B_{j}+\sigma_{-1}(B_{i}y_{j}+B_{j}y_{i})+\sigma_{-2}y_{i}y_{j},
\end{equation}
where we use the notations
\begin{equation}
\sigma=\frac{1}{\rho},\;\; \sigma_{0}=\frac{\rho_{0}}{\rho\tau},\;\;\; \sigma_{-2}=\frac{\rho_{-2}}{\rho\tau}.
\end{equation}
\end{proposition}
We can get another result from the Lemma (3.4) such that
\begin{theorem}
The  Cartan tensor $C^{ijk}$ of a Cartan space $\ell^n$ with $(\alpha,\beta)$-metric is given by
\begin{eqnarray}
C^{ijk}=\frac{-1}{2}[r_{-1}b^{i}b^{j}b^{k}+\varPi _{ijk}\lbrace\rho_{-1}a^{ij}b^{k}+\rho_{-2}a^{ij}Y^{k}+r_{-2}b^{i}b^{j}Y^{k} \\\nonumber + r_{-3}b^{i}Y^{j}Y^{k}\rbrace  + r_{-4}Y^{i}Y^{j}Y^{k}]
\end{eqnarray}
\end{theorem}
where the notation $\varPi_{ijk}$ means the cyclic symmetrization of the quantity in the brackets with respect to indices i, j, k.\\

We can deduce the other important geometric object fields for $\ell^{n}$ with $(\alpha,\beta)$-metric, for instance, $N_{ij}$,$H_{jk}^{i}$,$C_{i}^{jk}$ etc. without difficulty.

\section{Cartan spaces with infinite series of $(\alpha,\beta)$-metric}

In 2004 Lee and Park \cite{LP} introduced a r-th series $(\alpha,\beta)$-metric
\begin{equation}
L(\alpha,\beta)=\beta\sum_{k=0}^{r}(\frac{\alpha}{\beta})^{k},
\end{equation}
where they assume $\alpha <\beta$.
If $ r = 1$ then $ L =\alpha +\beta$ is a Randers metric. If $ r = 2$ then $ L = \alpha + \beta +\frac{\alpha^{2}}{\beta}$ is a combination of Randers metric and Kropina metric. If $ r = \infty$
then above metric is expressed as  
\begin{equation}
L(\alpha,\beta)=\frac{\beta^{2}}{\beta-\alpha}
\end{equation}

and the metric (5.2) named as infinite series $(\alpha, \beta)$-metric. This metric is very remarkable because it is the difference of Randers and Matsumoto metric.\\

In this section we consider two cases of Cartan Finsler spaces with special $(\alpha,\beta)$-metrics of deformed infinite series metric which are defined as

\begin{flushleft}
	\textbf{I} $ \;\;\; H(\alpha,\beta)=\frac{\alpha\beta^{2}}{\beta-\alpha} $ i.e. the product of infinite series and Riemannian metric.  \\
	\textbf{II} $ \;\; H(\alpha,\beta)=\frac{\beta^{3}}{\beta-\alpha} $ i.e. the product of infinite series and one-form metric.	
\end{flushleft}
	 
\subsection{Cartan space $\ell^{n}$ for $H(\alpha,\beta)=\frac{\alpha\beta^{2}}{\beta-\alpha}$}
In the first case, partial derivatives of the fundamental function $H(\alpha,\beta)$ lead us the followings:\\
\begin{center}
$H_{\alpha}=\frac{\beta^{3}}{(\beta-\alpha)^{2}}$,\;\;\; $H_{\beta}=\frac{\alpha\beta^{2}-2\alpha^{2}\beta}{(\beta-\alpha)^{2}}$\\
$H_{\alpha\alpha}=\frac{2\beta^{3}}{(\beta-\alpha)^{3}}$\;\;\;$H_{\alpha\beta}=\frac{\beta^{3}-3\alpha\beta^{2}}{(\beta-\alpha)^{3}}$\;\;\;$H_{\beta\beta}=\frac{2\alpha^{3}}{(\beta-\alpha)^{3}}$\\
$H_{\alpha\alpha\alpha}=\frac{6\beta^{3}}{(\beta-\alpha)^{4}}$\;\;\;$H_{\alpha\alpha\beta}=\frac{-6\alpha\beta^{2}}{(\beta-\alpha)^{4}}$\;\;\;$H_{\beta\beta\beta}=\frac{-6\alpha^{3}}{(\beta-\alpha)^{4}}$\;\;\;$H_{\alpha\beta\beta}=\frac{6\alpha^{2}\beta}{(\beta-\alpha)^{4}}$
\end{center}
Using equation (3.15) and (3.17) we have following invariants
\begin{equation}
\rho_{1}=\frac{\alpha\beta^{2}-2\alpha^{2}\beta}{2(\beta-\alpha)^{2}},\;\;\rho=\frac{\beta^{3}}{2(\beta-\alpha)^{2}},\;\;\rho_{0}=\frac{\alpha^{3}}{(\beta-\alpha)^{3}},\;\;\; 
\end{equation}
\begin{equation}
 \rho_{-1}=\frac{\beta^{3}-3\alpha\beta^{2}}{2\alpha(\beta-\alpha)^{3}}\;\;\; \rho_{-2}=\frac{\beta^{3}(3\alpha-\beta)}{2\alpha^{3}(\beta-\alpha)^{3}}\nonumber
\end{equation}
\begin{proposition}
The invariants $\rho$ never vanishes in a Cartan space $\ell^{n}$ equipped with deformed infinite series metric function $H(\alpha,\beta)=\frac{\alpha\beta^{2}}{\beta-\alpha}$ -metric on $\widetilde{T^{*}M}$. Converesely, we have $H_{\alpha}\ne 0$ on $\widetilde{T^{*}M}$.
\end{proposition}
Again using equation (3.20) and (3.23) we have following invariants
\begin{equation}
 r_{-1}=\frac{-3\alpha^{3}}{(\beta-\alpha)^{4}},\;\;\;\;\;\;\;\; r_{-2}=\frac{3\alpha\beta}{(\beta-\alpha)^{4}}, 
\end{equation}
\begin{equation}
r_{-3}=\frac{4\alpha\beta^{3}-\beta^{4}-9\alpha^{2}\beta^{2}}{2\alpha^{3}(\beta-\alpha)^{4}}\;\;\;r_{-4}=\frac{15\alpha^{2}\beta^{3}+3\beta^{5}-12\alpha\beta^{4}}{2\alpha^{5}(\beta-\alpha)^{4}}\nonumber
\end{equation}
\begin{proposition}
	The  invariants of Cartan tensor $ C^{ijk} $ in Cartan space $\ell^{n}$ which equipped with deformed infinite series metric function $H(\alpha,\beta)=\frac{\alpha\beta^{2}}{\beta-\alpha}$  is given by (5.2).
\end{proposition}
The invariants of equations (5.1) and (5.2) satisfies the following relations
\begin{equation}
\alpha^{2}\rho_{-1}+\beta\rho_{0}=\rho_{1},\;\;\;\alpha^{2}\rho_{-2}+\beta\rho_{-1}=0,
\end{equation}
\begin{equation}
\alpha^{2} r_{-2}+\beta r_{-1}=0,\;\;\;\alpha^{2} r_{-3}+\beta r_{-2}=-\rho_{-1},\nonumber
\end{equation}
\begin{equation}
\alpha^{2} r_{-4}+\beta r_{-3}=-2\rho_{-2},\nonumber
\end{equation}
\begin{theorem}
The Cartan space $\ell^{n}$ equipped with deformed infinite series metric function $H(\alpha,\beta)=\frac{\alpha\beta^{2}}{\beta-\alpha}$ has the invariants in equations (5.1) and (5.2) are satifies the relations in (5.3).
\end{theorem}
From equation (5.1) and (5.2), The fundamental tensor $g^{ij}(x,y)$ is of the form
\begin{equation}
g^{ij}(x,y)=\frac{\beta^{3}}{2(\beta-\alpha)^{2}} a^{ij}+\frac{\alpha^{3}}{(\beta-\alpha)^{3}} b^{i}b^{j}+\frac{\beta^{3}-3\alpha\beta^{2}}{2(\beta-\alpha)^{3}}(b^{i}Y^{j}+b^{j}Y^{i})
\end{equation}
\begin{equation}
+\frac{\beta^{3}(3\alpha-\beta)}{2\alpha^{3}(\beta-\alpha)^{3}}Y^{i}Y^{j}\nonumber
\end{equation}
\begin{corollary}
The fundamental tensor $g^{ij}(x, y)$ of the space $\ell^{n}$ endoewd with the metric function $H(\alpha,\beta)=\frac{\alpha\beta^{2}}{\beta-\alpha}$ is given by the equation (5.4).\\
\end{corollary}
Converesely we obtain
\begin{theorem}
The Cartan space with $(\alpha,\beta)$-metric which have the invariants such that (5.1)and (5.2) is the spaces $\ell^{n}$ with the fundamental function $H(\alpha,\beta)=\frac{\alpha\beta^{2}}{\beta-\alpha}$, i.e. deformed infinite series metric. 
\end{theorem}
\subsection{Cartan space $\ell^{n}$ for $H(\alpha,\beta)=\frac{\beta^{3}}{\beta-\alpha}$}

In the second case, partial derivatives of the fundamental function $H(\alpha,\beta)$ lead us the followings:\\
\begin{center}
$H_{\alpha}=\frac{\beta^{3}}{(\beta-\alpha)^{2}}$,\;\;\; $H_{\beta}=\frac{2\beta^{3}-3\alpha\beta^{2}}{(\beta-\alpha)^{2}}$\\
$H_{\alpha\alpha}=\frac{2\beta^{3}}{(\beta-\alpha)^{3}}$\;\;\;$H_{\alpha\beta}=\frac{\beta^{3}-3\alpha\beta^{2}}{(\beta-\alpha)^{3}}$\;\;\;$H_{\beta\beta}=\frac{2\beta^{3}-6\alpha\beta^{3}+6\alpha^{2}\beta}{(\beta-\alpha)^{3}}$\\
$H_{\alpha\alpha\alpha}=\frac{6\beta^{3}}{(\beta-\alpha)^{4}}$\;\;\;$H_{\alpha\alpha\beta}=\frac{-6\alpha\beta^{2}}{(\beta-\alpha)^{4}}$\;\;\;$H_{\beta\beta\beta}=\frac{-6\alpha^{3}}{(\beta-\alpha)^{4}}$\;\;\;$H_{\alpha\beta\beta}=\frac{6\alpha^{2}\beta}{(\beta-\alpha)^{4}}$
\end{center}
Using equation (3.15) and (3.17) we have following invariants
\begin{equation}
 \rho_{1}=\frac{2\beta^{3}-3\alpha\beta^{3}}{2(\beta -\alpha)^{2}},\;\;\rho=\frac{\beta^{3}}{2(\beta-\alpha)^{2}},\;\;\rho_{0}=\frac{\beta^{3}-3\alpha \beta^{2}+3\alpha^{2}\beta}{(\beta-\alpha)^{3}},\;\;\; 
\end{equation}
\begin{equation}
 \rho_{-1}=\frac{\beta^{3}-3\alpha\beta^{2}}{2\alpha(\beta-\alpha)^{3}}\;\;\; \rho_{-2}=\frac{\beta^{3}(3\alpha-\beta)}{2\alpha^{3}(\beta-\alpha)^{3}}\nonumber
\end{equation}
\begin{proposition}
	The invariants $\rho$ never vanishes in a Cartan space $\ell^{n}$ equipped with deformed infinite series metric function $H(\alpha,\beta)=\frac{\beta^{3}}{\beta-\alpha}$ -metric on $\widetilde{T^{*}M}$. Converesely, we have $H_{\alpha}\ne 0$ on $\widetilde{T^{*}M}$.
\end{proposition}
Again using equations (3.20) and (3.23) we have following invariants
\begin{equation}
 r_{-1}=\frac{-3\alpha^{3}}{(\beta-\alpha)^{4}},\;\;\;\;\;\;\;\; r_{-2}=\frac{3\alpha\beta}{(\beta-\alpha)^{4}}, 
\end{equation}
\begin{equation}
r_{-3}=\frac{4\alpha\beta^{3}-\beta^{4}-9\alpha^{2}\beta^{2}}{2\alpha^{3}(\beta-\alpha)^{4}}\;\;\;r_{-4}=\frac{15\alpha^{2}\beta^{3}+3\beta^{5}-12\alpha\beta^{4}}{2\alpha^{5}(\beta-\alpha)^{4}}\nonumber
\end{equation}
\begin{proposition}
	The  invariants of Cartan tensor $ C^{ijk} $ in Cartan space $\ell^{n}$ which equipped with deformed infinite series metric function $H(\alpha,\beta)=\frac{\beta^{3}}{\beta-\alpha}$ is given by (5.6).
\end{proposition}
The invariants of equations (5.5) and (5.6) satisfies the following relations
\begin{equation}
\alpha^{2}\rho_{-1}+\beta\rho_{0}=\rho_{1},\;\;\;\alpha^{2}\rho_{-2}+\beta\rho_{-1}=0,
\end{equation}
\begin{equation}
\alpha^{2} r_{-2}+\beta r_{-1}=0,\;\;\;\alpha^{2} r_{-3}+\beta r_{-2}=-\rho_{-1},\nonumber
\end{equation}
\begin{equation}
\alpha^{2} r_{-4}+\beta r_{-3}=-2\rho_{-2},\nonumber
\end{equation}
\begin{theorem}
The Cartan space $\ell^{n}$ equipped with deformed infinite series metric function $H(\alpha,\beta)=\frac{\beta^{3}}{\beta-\alpha}$ has the invariants in equations (5.5) and (5.6) are satifies the relations in (5.7).
\end{theorem}
From equation (5.5) and (5.6), The fundamental tensor $g^{ij}(x,y)$ is of the form
\begin{equation}
g^{ij}(x,y)=\frac{\beta^{3}}{2(\beta-\alpha)^{2}} a^{ij}+\frac{\beta^{3}-3\alpha \beta^{2}+3\alpha^{2}\beta}{(\beta-\alpha)^{3}} b^{i}b^{j}
\end{equation}
\begin{equation}
+\frac{\beta^{3}-3\alpha\beta^{2}}{2(\beta-\alpha)^{3}}(b^{i}Y^{j}+b^{j}Y^{i})+\frac{\beta^{3}(3\alpha-\beta)}{2\alpha^{3}(\beta-\alpha)^{3}}Y^{i}Y^{j}\nonumber
\end{equation}
\begin{corollary}
The fundamental tensor $g^{ij}(x,y)$ of the space $\ell^{n}$ endoewd with the metric function $H(\alpha,\beta)=\frac{\beta^{3}}{\beta-\alpha}$ is given by the equation (5.8).
\end{corollary}
Converesely we obtain
\begin{theorem}
The Cartan space with $(\alpha,\beta)$-metric which have the invariants such that (5.5)and (5.6) is the spaces $\ell^{n}$ with the fundamental function $H(\alpha,\beta)=\frac{\beta^{3}}{\beta-\alpha}$,i.e. deformed infinite series metric. 
\end{theorem}

\section{Conclusions}

In this work, we consider the  infinite series $(\alpha,\beta)$-metric, Riemannian metric and 1-form metric we determine relations with the "invarints" which characterize the special classes in Cartan Finsler frames . But, in Finsler geometry, there are many($\alpha$,$\beta$)-metrics, in future work we can determine the frames for them also.

\pagebreak


\begin{thebibliography}{20}

\bibitem{Be}\textbf{Berwald, L.:} $ \ddot{U} $ber Finslersche und Cartansche Geometric I,II,III,IV, \textit{Mathematica Cluj}, 16 \textbf{(1940)}; \textit{Compositio Mathematica}, 7 \textbf{(1939)},  141-176; \textit{Ann. Math.}, 42\textbf{(1941)}, 84-112; 48\textbf{(1947)}, 755-781.

\bibitem{Br}\textbf{Brickell, F.:} A relation between Finsler and Cartan structures, \textit{Tensor, N. S.}, 25 \textbf{(1972)}, 360-364.

\bibitem{Car}\textbf{Cartan, E.:} Less spaces metriwues fondes sur la notion d$ ^{'} $aire, \textit{Actuarites Sci. Ind.}, no. 72, Hermann,Paris \textbf{1933}.

\bibitem{Iga1} \textbf{Igarachi, T.:} Lie derivatives in Cartan spaces, \textit{Tensor, N. S.}, 52 (\textbf{1993)}, 80-89.

\bibitem{Iga2}\textbf{Igarachi, T.:} $(\alpha,\beta)$ - Metrics in Cartan Spaces, \textit{Tensor, N. S.}, 55 \textbf{(1994)}, 74-83.

\bibitem{Kro}\textbf{Kropina, V. K.:} On the projective two-dimensional Finsler spaces with a special metric,(in Russian),\textit{Trudy Seminara po Vektornomu i Tenzornomu Analizu}, 11 \textbf{(1961)}, 277-292.

\bibitem{M1} \textbf{Matsumoto, M.:} Foundation of Finsler geometry and special Finsler spaces, Kaiseisha Press, \textbf{1986}.

\bibitem{M2} \textbf{Matsumoto, M.:} A slop of a mountain is a Finsler surface with respect to time measure, \textit{J. Math. Kyoto Univ.}, 29 \textbf{(1989)}, 17-25.

\bibitem{M3} \textbf{Matsumoto, M.:} Theory of Finsler spaces with $(\alpha,\beta)$-metric, \textit{Reports on Mathematical Physics}, 31\textbf{ (1992)}, 43-83.

\bibitem{Mi1}\textbf{Miron, R.:} Cartan spaces in a new point of view by considering them as duals of Finsler spaces, \textit{Tensor, N. S.}, 46 \textbf{(1987)}, 329-334.

\bibitem{Mi2}\textbf{Miron, R.:} The geometry of Cartan spaces, \textit{Progress of Math.}, 22 \textbf{(1988)}, 1-38.

\bibitem{Ra}\textbf{Randes, G.:} On an asymmetric metric in the four-space of general relativity, \textit{Physical Rev.}(2), 59 \textbf{(1941)}, 195-199.

\bibitem{Ru}\textbf{Rund, H.:} The Hamilton-Jacobi theory in the Calculus of Varitions,\textit{D.Van Nostrond Co.}, London, \textbf{1966}.

\bibitem{LP}\textbf{Lee, I. Y. and Park, H. S.:} Finsler spaces with infinite series $(\alpha,\beta)$- metric, {J.Korean Math.Society}, 41, (3), \textbf{(2004)}, 567-589.

\end{thebibliography}
\end{document}